\documentclass[11pt,reqno,oneside]{amsart}
\usepackage{latexsym}
\usepackage{graphics}
\author{F. Boniver and P. Mathonet}
\title{IFFT-equivariant quantizations}
\newtheorem{lem}{Lemma}
\newtheorem{thm}[lem]{Theorem}
\newtheorem{prop}[lem]{Proposition}
\newtheorem{cor}[lem]{Corollary}
\theoremstyle{remark}

\theoremstyle{definition}

\newcommand{\R}{\mathbb{R}}
\newcommand{\C}{\mathbb{C}}
\newcommand{\K}{\mathbb{K}}
\newcommand{\N}{\mathbb{N}}
\renewcommand{\O}{\mathfrak{O}^{(n)}}
\renewcommand{\Sp}{\mathfrak{S}^{(n)}}
\renewcommand{\L}{L}
\newcommand{\D}{\mathcal{D}}
\renewcommand{\S}[1][\delta]{\mathcal{S}_{{#1}}}
\newcommand{\F}{\mathcal{F}}
\newcommand{\g}{\mathfrak{g}}
\newcommand{\h}{\mathfrak{h}}
\renewcommand{\H}{\mathfrak{H}}
\newcommand{\euler}{\mathcal{E}}
\newcommand{\ad}{\mathit{ad}}
\newcommand{\Lop}[1][]{\L^{\mathit{op}{#1}}}
\newcommand{\Lt}{\L^t}
\newcommand{\del}[1]{\frac{\partial}{\partial x^{{#1}}}}
\newcommand{\casi}[2]{\mathcal{C}_{{#2}}^{#1}}
\newcommand{\casinil}{N_\mathcal{C}}
\newcommand{\tr}{\mathrm{tr}}
\newcommand{\aff}{\mathit{Aff}}

\newcommand{\restr}[1]{\mbox{}_{\vert {#1}}}
\newcommand{\Vect}{\mathrm{Vect}}
\newcommand{\alg}[2]{\mathit{{#1}}({#2})}
\newcommand{\rhot}[1]{\rho^{{#1}}}

\newcommand{\tree}[1][\gamma]{\mathcal{T}_{{#1}}}
\newcommand{\treebig}{\widetilde{\tree[]}\,}
\newcommand{\smooth}{smooth\ }
\begin{document}
\begin{abstract}
The existence and uniqueness of quantizations that are equivariant with respect to conformal and projective Lie algebras of vector fields were recently obtained by Duval, Lecomte and Ovsienko. In order to do so, they computed spectra of some Casimir operators. We give an explicit formula for those spectra in the general framework of \emph{IFFT}-algebras classified by Kobayashi and Nagano. We also define \emph{tree-like} subsets of eigenspaces of those operators in which eigenvalues can be compared to show the existence of IFFT-equivariant quantizations.  We apply our results to prove existence and uniqueness of quantizations that are equivariant with respect to the infinitesimal action of the symplectic (resp. pseudo-orhogonal) group on the corresponding  Grassmann manifold of maximal isotropic subspaces.
\end{abstract}
    \maketitle
\noindent
Math. Classification (AMS 2000): 17B66, 22E46, 81R05. 
\section{Introduction}
 The word ``quantization'' carries several different meanings, both in physics and mathematics.
    One approach --- see for instance~\cite{wo} --- is to consider a quantization procedure as a linear bijection 
    from the space of \emph{symbols} $\mbox{Pol}(T^{\star}M)$ of \smooth functions on the cotangent 
    bundle of a manifold $M$ that are polynomial along the fibres 
    to the space $\mathcal{D}_{\frac{1}{2}}(M)$ 
    of linear differential operators acting on half-densities.
    It is known that these spaces cannot be canonically identified. 
    In other words, there does not exist a preferred quantization 
    procedure.
    
    The concept of \emph{equivariant quantization} was introduced and developped in 
    \cite{lo,lecras} and \cite{dlo} by C. Duval, P.B.A. Lecomte and V. Ovsienko. 
    These recent works  take care of the 
    symmetries of the classical situation to quantize.
    
    If $G$ is a group acting on 
    the manifold $M$, a $G$-equivariant quantization consists in an  isomorphism of representations of $G$ between 
    the spaces of symbols and of differential operators.
    Obviously, such an identification does not exist for all groups $G$ acting 
    on $M$ : for instance those spaces are not equivalent as $\mathit{Diff}(M)$-modules.
     At the infinitesimal level, if $G$ is a Lie group, its  action gives rise to a Lie subalgebra $\g$ of vector fields over $M$ and one is lead to building a $\g$-equivariant linear bijection.
    Lecomte and Ovsienko examined the case of a projective structure on a manifold of 
    dimension $n$, with $G = SL(n+1,R)$ and then, together with Duval, the case of 
    the group $G = SO(p+1,q+1)$  on a manifold of dimension $p+q$.  That latter group defines conformal transformations with respect to a pseudo-Riemannian metric.
    
    In those works, the authors consider the more general modules $\mathcal{D}_{\lambda,\mu}$ of differential operators transforming $\lambda$-densities into $\mu$-densities.  These parameters give rise to the \emph{shift} value $\delta=\mu-\lambda$ and to the special case $\delta=0$, which can be specialized to the original problem.  They obtain existence and uniqueness (up to normalization) results for  a 
    quantization procedure in both projective and conformal cases, provided the shift value does not belong to a \emph{critical} set.  Furthermore, they show that this set never contains zero.

    In suitable charts, the subalgebras mentioned up to now are realized by polynomial vector fields.  Furthermore, they  share the 
    property of being maximal in the algebra of polynomial 
    vector fields. 
    
    In \cite{bm}, we investigated this maximality property and showed 
    that the finite dimensional, graded and maximal subalgebras of 
    the Lie algebra of polynomial vector fields over a Euclidean vector space correspond 
    to the list of so called ``Irreducible Filtered Lie algebras of 
    Finite Type''(or IFFT-algebras), classified by S. Kobayashi and T. Nagano in~\cite{koba1}.
    
    Our concern in this paper is to deal with the natural next question~:
    \emph{``Do there exist (unique) equivariant quantizations with respect to the IFFT-algebras ?'' }
    
    The original construction of the conformally equivariant 
    quantization (see~\cite{dlo}) involves the computation
    of the spectrum of the Casimir operator
    of $so(p+1,q+1)$ acting on the space of symbols.
    Obstructions to its existence show up as equalities among 
    some eigenvalues of that operator. It was also shown how the 
relevant eigenvalues that should be compared are associated to 
\emph{tree-like subsets}\  of eigenspaces.

     Section~\ref{sec.computcas} of the present article is devoted to this computation.
    We obtain, for a wide range of IFFT-algebras, a formula 
where the eigenvalues are expressed 
    in terms of the dimension of the manifold and of the highest 
    weights of some finite dimensional representations of the semisimple part of the linear isotropy algebra of $\g$~(see~\cite{koba3}).
    
    In Section~\ref{sec.buildeqq}, we propose a general definition for the above mentioned tree-like subsets.  A few elementary properties of these subspaces allow us to reformulate the existence theorem for equivariant quantizations in the framework of IFFT-algebras.
    
    We eventually apply those results in Section~\ref{sec.examples}.
    The Lie algebras of fundamental vector fields associated to the action
    of the symplectic (resp. pseudo-orthogonal) group on the Lagrangian
    (resp. pseudo-orthogonal) Grassmann manifold are indeed IFFT.
    We prove existence and uniqueness results for equivariant quantizations with respect to both 
of those algebras. Once more, these results hold outside of a critical set of values of the shift.
  We furthermore prove that this set never contains zero.

\section{Definitions and notations}
The basic geometric objects that are involved in the present work may be defined in a standard way over an arbitrary manifold.  Those definitions were for instance  given in~\cite{dlo}.  For the sake of completeness and in order to fix some notations, we shall briefly recall them.  It will be sufficient for  our computations to explicit those over  vector spaces. 

Throughout this Section, let $V$ denote a $d$-dimensional vector space over  $\K=\R$ or $\C$. 

Whenever $E$ is a vector bundle over $V$, the space of \emph{sections} of $E$, which we shall write $\Gamma(E)$, is taken to be the space of $C_\infty$ sections if $\K=\R$ or the space of holomorphic sections if $\K=\C$.  

\subsection{Differential operators acting on densities}
Strictly speaking, a \emph{density of weight $\lambda\in\R$} or \emph{$\lambda$-density} over $V$ is a function
\[
\phi : \wedge^d V\setminus\{ 0\} \rightarrow \K 
\]
that is homogeneous of degree $\lambda$, i.e. such that 
\[
\phi(s \omega)=\vert s\vert^\lambda \phi(\omega),\quad \forall s\in\K\setminus\{0\}, \forall \omega\in\wedge^d V\setminus\{0\}.
\]
The one-dimensional vector space of $\lambda$-densities is then denoted by $\Delta^\lambda V$.

However, we shall also name density a section of $V\times\Delta^\lambda V$ viewed as a trivial vector bundle over $V$.
The space $\F_\lambda$ of such sections carries a natural structure of module over the Lie algebra $\Vect(V)$ of vector fields over $V$.  The Lie derivative is given by
\[
\L^\lambda_X \phi = X.\phi + \lambda \,\tr(DX) \phi,\quad \forall X\in\Vect(V),\forall \phi\in\F_\lambda,
\]
where $DX$ denotes the Jacobian matrix of $X$.

Let now $\D_{\lambda, \mu}$ be the space of linear differential operators from $\F_\lambda$ to $\F_\mu$.  It is a module over $\Vect(V)$ as well. The Lie derivative $\Lop[\lambda, \mu] D$ of a differential operator $D$ is its natural derivative as a linear map from $\F_\lambda$ to $\F_\mu$:
\[
\Lop[\lambda,\mu] D=\L^\mu\circ D - D\circ \L^\lambda.
\]
In order to keep light notations, we shall simply write $\Lop$ for $\Lop[\lambda, \mu]$ unless that leads to confusion.  

To the module $\D_{\lambda,\mu}$ is associated the \emph{shift}\ value  $\delta=\mu-\lambda$.

The space $\D_{\lambda, \mu}$ is the filtered union $\bigcup_{k\in\N} \D_{\lambda, \mu}^k$ of the spaces of differential operators of order at most $k$. Clearly, the module structure defined by $\Lop$ restricts itself to these.
\subsection{Symbols}
Let $\vee^k V$ denote the $k$-th symmetric power of $V$.  The principal symbol of $D\in\D_{\lambda,\mu}^k$ is a section  of
\[
V \times (\vee^k V\otimes \Delta^\delta V),
\]
viewed  as a trivial bundle over $V$.  
The space $\S^k$ of such sections is called \emph{symbol space of degree $k$}.
We also consider the whole symbol space
\[
\S=\bigoplus_{k\geq 0} \S^k.
\]

The module $(\S, \Lt)$, where $\Lt$ is the natural Lie derivative on tensor fields valued in $\delta$-densities, is then the graded module associated to $(\D_{\lambda,\mu},\Lop)$.
\subsection{Basic identifications}\label{subsec.bid}
The space $\vee^k V$ is naturally identified to the space of homogeneous polynomials  on $V^*$ of degree $k$.

Let $c\in\F_\delta$ and $T\in\vee^k V$.  The tensor field $c T\in\S^k$ may be described as a function on $T^* V\cong V\times V^*$:
\[
cT:(x,\xi)\in V\times V^*\mapsto c(x)T(\xi).
\]
The latter expression is a polynomial of degree $k$ in $\xi$.  
Throughout this paper, we shall consider $\xi$ as the indeterminate variable to write such tensor fields as ``polynomial'' functions.

With these notations, the Lie derivative of a tensor field $P\in\S$ can be explicited as follows:
\begin{equation}\label{eq.lt}
\Lt_X P = \sum_i X^i \del{i} P-\sum_{i,j} \del{j}X^i \xi_i \frac{\partial}{\partial \xi_j} P+\delta \sum_i \del{i} X^i P.
\end{equation}

Moreover, we will use this notation for differential operators as well.  Any $D\in\D_{\lambda,\mu}$
 may be written
\[
f\in\F_\lambda \mapsto \sum_{\vert\alpha\vert\leq k} c_\alpha D^\alpha f\in\F_\mu,
\]
where $\alpha$ is a multi-index, $D^\alpha$ stands for $(\del{1})^{\alpha_1}\cdots(\del{d})^{\alpha_d}$ and $c_\alpha\in\F_\delta$.

The \emph{principal symbol} of $D$ is
\begin{equation}
  \label{eq.prsymb}
\sigma(D)=  \sum_{\vert\alpha\vert = k} c_\alpha \xi^\alpha
\end{equation}
and the \emph{total symbol map}, which is also known as  the \emph{normal ordering}, 
\[
\sigma^{\lambda,\mu}_{\aff}:\D_{\lambda,\mu}\to \S: D\mapsto \sum_{\vert\alpha\vert\leq k} c_\alpha \xi^\alpha
\]
is a vector space isomorphism. 
The index $\aff$ recalls that $\sigma^{\lambda,\mu}_{\aff}$ is moreover equivariant with respect to the \emph{affine} Lie subalgebra $\aff$, which is made of constant and linear vector fields (cf. Proposition \ref{prop.lieder}).

We can also endow $\S$ with the module structure that makes $\sigma^{\lambda,\mu}_{\mathit{\aff}}$ a module isomorphism.  This is done by considering the representation 
\[
\sigma^{\lambda,\mu}_\aff \circ \Lop[\lambda,\mu] \circ (\sigma^{\lambda,\mu}_\aff)^{-1},
\]
which we still denote $\Lop[\lambda,\mu]$ or simply $\Lop$.
The comparison of spaces of differential operators and tensor fields as modules over a given subalgebra of vector fields becomes the comparison of the modules $(\S,\Lt)$ and $(\S,\Lop)$, provided one keeps in mind that two parameters, namely $\lambda$ and $\mu$, are attached to the second one.

\subsection{Equivariance algebras}
Let us now present the subalgebras that define the module structures of  $(\S,\Lt)$ and $(\S,\Lop)$.
Those are the \emph{Irreducible Filtered Finite-dimensional Transitive} Lie algebras, listed by Kobayashi and Nagano in~\cite{koba1}.  One of the main results is that any  IFFT-algebra $\g$ carries a canonical \emph{$3$-grading}:
\[
\g=\g_{-1}\oplus \g_0 \oplus \g_1.
\] 
We shall use the following properties of an IFFT-algebra $\g$ presented with its decomposition.
\begin{itemize}
\item $\g$ is simple.
\item $\g_0$ is reductive and is the sum of its semisimple part $\h_0$ and its one dimensional center spanned by the \emph{Euler} element $\euler$.
\item $\g_p$ is the eigenspace of eigenvalue $p$ of $\ad(\euler)$.
\item The subspaces $\g_{-1}\oplus\g_1$ and $\g_0$ are orthogonal with respect to the Killing form $\beta$ of $\g$.
\end{itemize}

Moreover, there exists a standard imbedding of $\g$ in the Lie algebra of polynomial vector fields  over $\g_{-1}$ (see~\cite{koba3}).  The grading of $\g$ then corresponds to the grading of polynomial vector fields with respect to the degree of the coefficients. In other words, $\g_{-1}$ is identified to the space of constant vector fields, while $\g_0$ and $\g_1$ correspond to some linear and quadratic vector fields respectively.

In~\cite{bm}, we proved that this imbedding is maximal if $\K=\C$ or if $\K=\R$ and $\g_{-1}$ admits no complex structure as a module over $\g_0$.

In the present paper, we shall compare the modules $(\S,\Lt)$ and $(\S,\Lop)$ over the base space $V=\g_{-1}$.

\subsection{Equivariant quantizations and symbol maps}
A \emph{($\g$-)equivariant symbol map} is a $\g$-module isomorphism 
\[
\sigma_\g : \D_{\lambda, \mu}\to \S
\]
that induces the identity on the associated graded module.   Explicitly, this latter requirement means that
\[
D\in \D_{\lambda, \mu}^k\implies \sigma_\g(D)-\sigma(D)\in \bigoplus_{l<k}\S^l,
\]
where $\sigma$ associates to a differential operator its principal symbol (see~(\ref{eq.prsymb})).
With respect to the identification described in Subsection~\ref{subsec.bid}, it is equivalent to obtain a map with similar properties between the modules $(\S,\Lop)$ and $(\S,\Lt)$.
The inverse map of such an application is named \emph{($\g$-)equivariant quantization}.

\subsection{Casimir operators}
Let now $(V,\rho)$ denote  a module over a semisimple Lie algebra $\g$ and $\beta$ a nondegenerate invariant bilinear form over $\g$.

Let $(e_i)$ $(i=1,\ldots, p)$ be a basis of $\g$ and $(e_i^*)$ $(i=1,\ldots, p)$ its dual basis with respect to $\beta$, which is defined by the conditions
\[
\beta(e_i, e_j^*)=\delta_{ij}.
\]
As it is usual, we name \emph{Casimir operator} associated to $\g, \rho$ and $\beta$,  the following endomorphism of $V$:
\[
\casi{\rho}{\g, \beta}=\sum_{i=1}^p \rho(e_i)\rho(e_i^*).
\]
Once the bilinear form $\beta$ is fixed, it does not depend on the choice of the basis $(e_i)$.
%%%%%%%%%%%%%%%%%%%%%%%%%%%%%%%%%%%%%%%%%
\section{Computing Casimir operators}\label{sec.computcas}%%%%
%%%%%%%%%%%%%%%%%%%%%%%%%%%%%%%%%%%%%%%%%
From now on, we shall assume that 
$\g$
is the standard maximal imbedding of an IFFT-algebra and $\beta$ its Killing form. We set
$d=\dim(\g_{-1})$ and denote by $\beta_0$ the Killing form of $\h_0$.  
\subsection{Some preliminaries}
The following result follows from a direct computation.
\begin{prop}\label{prop.lieder}
  The modules $(\S,\Lop)$ and $(\S,\Lt)$ coincide over the Lie algebra $\aff$.
\end{prop}

We shall now describe suitable bases of $\g$ in  order to simplify the computation of the Casimir operator.
\begin{prop}\label{prop.dualbases}%%%%%%%%%%%%%%%%%   Choix de bases duales appropriees
  Let $(e_i)$ $(i=1,\ldots, d)$ denote a basis of $\g_{-1}$ and $(h_j)$ $(j=1,\ldots, \dim(\h_0))$ a basis of $\h_0$.  There exist unique bases $(\epsilon^i)$ and $(h_j^*)$ of $\g^1$ and $\h_0$ respectively such that the bases $(e_i,\euler,h_j,\epsilon^i)$ and $(\epsilon^i,\frac{1}{2 d} \euler, h_j^*, e_i)$ of $\g$ are dual to each other with respect  to $\beta$.

Moreover, 
\begin{equation}\label{eq.eeps}
\sum_i [e_i,\epsilon^i]=-\frac{1}{2}\euler.
\end{equation}
\end{prop}
\begin{proof}
  The restriction of $\beta$ to both spaces $\g_{-1}\oplus\g_1$ and $\g_0$ is nondegenerate.
  Hence the existence and uniqueness of the $(\epsilon^i)$ basis in $\g_1$.
  
  But $\h_0$ and $\K \euler$ are orthogonal to each other too.  Indeed,
  any $h\in \h_0$ can be written $\sum_{r=1}^R [h_r, h'_r]$ for
  some $R\in\N$ and $h_r, h'_r\in\h_0$ since $\h_0$ is semisimple. Then,
  for all $k,l \in \K$ and $h\in\h_0$,
\begin{eqnarray}
  \beta(\euler, k\euler + l h)
  &=& k \beta(\euler, \euler) - l \sum_r \beta([h_r, \euler], h'_r) \nonumber \\
  &=& k \tr(\ad(\euler)^2) \nonumber \\
  &=& 2k\, d.\label{eq.euleul}
\end{eqnarray}
  Hence the existence and uniqueness of the $(h_j^*)$ basis and the value of the dual element of $\euler$.

Moreover, for all $k_1,k_2\in\g_0$,
\begin{eqnarray*}
  \beta([\sum_i [e_i,\epsilon^i], k_1],k_2) &=& \sum_i \beta(\epsilon^i,[[k_1,k_2],e_i])\\
  &=& \tr(\ad([k_1,k_2])\restr{\g_{-1}}) \\
  &=& \tr([\ad(k_1)\restr{\g_{-1}},\ad(k_2)\restr{\g_{-1}}]) \\
  &=& 0.
\end{eqnarray*}
Therefore, $\sum_i [e_i,\epsilon^i]$ is in the center of $\g_0$ and can thus be written $r\euler$ with $r\in \K$.
But, because of~(\ref{eq.euleul}),
\begin{eqnarray*}
  2r\, d &=& \beta(r\euler,\euler) \\
  &=& \beta(\sum_i [e_i, \epsilon^i],\euler)\\
  &=&-\sum_i \beta(\epsilon^i, e_i)\\
  &=&-d.
\end{eqnarray*}
Hence the conclusion.
\end{proof}
%%%%%%%%%%%%%%%%%%%%%%%%%%%% FIN Choix de bases duales appropriees
%%
\subsection{The cocycle $\gamma$}
According to Proposition~\ref{prop.lieder}, the ``only'' difficulties to build a $\g$-equivariant quantization stand in $\g_1$.  They are best seen in the difference of the Casimir operators on differential operators and symbols.  We shall denote by $\gamma$ the map
\[
\g\rightarrow \alg{gl}{\S}: X\mapsto \Lop_X-\Lt_X.
\]
This map is a Chevalley-Eilenberg cocycle with values in $\alg{gl}{\S}$, the latter representation being endowed with the Lie derivative defined by
\[
L_X :\alg{gl}{\S}\rightarrow \alg{gl}{\S}: T\mapsto \Lop_X \circ T- T\circ \Lt_X.
\]
Using the notations of Proposition~\ref{prop.dualbases}, we state the following observation.
\begin{prop}\label{prop.casi1}
  The Casimir operators of $\D_{\lambda,\mu}$ and $\S$ are related by the formula
\[
\casi{\Lop}{\g, \beta} = \casi{\Lt}{\g, \beta} +2 \sum_i \gamma({\epsilon^i})\circ \Lt_{e_i}.
\]
\end{prop}
\begin{proof}
  The Casimir operator $\casi{\Lop}{\g,\beta}$ can be rewritten as follows:
  \begin{eqnarray*}
    \casi{\Lop}{\g,\beta} &=& 
    \sum_i(\Lop_{e_i}\circ\Lop_{\epsilon^i}+\Lop_{\epsilon^i}\circ\Lop_{e_i}) 
    +\frac{1}{2d}(\Lop_{\euler})^2+\sum_j \Lop_{h_j}\circ\Lop_{h_j^*}\\
    &=& 2\sum_i \Lop_{\epsilon^i}\circ \Lop_{e_i} + \Lop_{\sum_i [e_i, \epsilon^i]}
    +\frac{1}{2d}(\Lop_{\euler})^2+\sum_j \Lop_{h_j}\circ\Lop_{h_j^*}.
  \end{eqnarray*}
  The conclusion is then a direct consequence of the vanishing of $\gamma$ on $\g_{-1}\oplus \g_0$.
\end{proof}
\subsection{Diagonalizing $\casi{\Lt}{\g,\beta}$}\label{subsec.diag}
We shall now specify the Casimir operator associated to the natural Lie derivative of tensor fields.

The Jacobian correspondence that associates to any linear vector field $X$ its coefficient matrix $\del{}X$ is an antihomomorphism from the Lie algebra of linear vector fields to $\alg{gl}{\g_{-1}}$.  Let us denote by $\H_0$ the subalgebra of $\alg{gl}{\g_{-1}}$ isomorphic to $\h_0$  through the opposite of the Jacobian and by $\rhot{k}$ the representation of $\h_0$ on the fiber $\vee^k \g_{-1}\otimes \Delta^\delta \g_{-1}$ of $\S^k$ which is induced by the restriction to $\H_0$ of the natural representation of $\alg{gl}{\g_{-1}}$ on this space.   Note that the latter does not depend on the value of $\delta$ since $\H_0$ is made of traceless matrices.  It is thus isomorphic to $\vee^k \g_{-1}$ endowed with the (restriction of) the standard representation of $\alg{gl}{\g_{-1}}$.
\begin{prop}\label{prop.casi2}
  Let $P\in\S^k$.  Then 
  \begin{equation}\label{eq.casit1}
    \casi{\Lt}{\g,\beta}P=\frac{1}{2d}(d\delta -k)(d(\delta -1)-k)P+\casi{\rhot{k}}{{\h_0},\beta\restr{\h_0}}P,
  \end{equation}
  where the last term is computed pointwise.
\end{prop}
\begin{proof}
  As in Proposition~\ref{prop.casi1}, we write
  \[
  \casi{\Lt}{\g,\beta}= \sum_i \Lt_{\epsilon^i}\circ\Lt_{e_i} -\frac{1}{2}\Lt_{\euler} +\frac{1}{2d}(\Lt_{\euler})^2 +\sum_j \Lt_{h_j}\circ\Lt_{h_j^*}.
  \]
%  Now, assume that the right hand side above has been completely developped and rewritten
%  \begin{equation}\label{eq.casiterms}
%    \sum_{(\alpha)} c_\alpha D_\alpha(P),
%  \end{equation}
%  with each $c_\alpha$ a \smooth function and $D_\alpha$ a linear differential operator in both horizontal and vertical variables, with constant coefficients.   
Being invariant under $\g$ and in particular $\g_{-1}$, $\casi{\Lt}{\g,\beta}$ is a differential operator with constant coefficients.  In order to compute it,  we thus only need to sum the constant terms in the right hand side of the last formula.
  In view of the explicit Formula~(\ref{eq.lt}) of $\Lt$, it is clear that the Lie derivatives with respect to a quadratic vector field do not contribute to such  terms.  Furthermore, it is easily checked that 
  \begin{equation}\label{eq.lieeuler}
  \Lt_{\euler}P=\sum_i x^i\del{i}P+(d\delta-k)P
  \end{equation}
  and that
  \[
  \Lt_{h_j}P=\sum_i (h_j)^i \del{i}P + \rhot{k}(h_j)P. 
  \]
  Hence the proof.
\end{proof}
The last result shows that, in order to compute $\casi{\Lt}{\g,\beta}$, it is sufficient to diagonalize a Casimir operator built upon a finite-dimensional representation.  We refer to~\cite{bou1} for results about semisimple Lie algebras.

Let us now introduce a few more notations needed to state the main theorem.
From now on to the end of this Section, for each vector space (resp. Lie algebra) $E$, we shall denote by $E^\C$ the complexified vector space (resp. Lie algebra) $E\otimes_\R \C$.  We shall set 
\[
\widetilde{E}=\left\{
\begin{array}{l}
E \text{ if the base field $\K$ is $\C$,}\\
E^\C \text{ if $\K=\R$.}
\end{array}
\right.
\]
Furthermore, we fix a Cartan subalgebra $\mathfrak{C}$ in $\widetilde{\h_0}$, a root system $\Lambda$,  a simple root system $\Lambda_S$.  Finally, let us denote by $S$ the sum of the positive roots and by $(\cdot,\cdot)$ the scalar product induced by the extension of $\beta_0$ to $\widetilde{\h_0}$ on
\[
\mathfrak{C}_r=\{h\in\mathfrak{C}:\lambda(h)\in\R,\forall \lambda\in\Lambda\}.
\]

If $E$ is an irreducible  module over $\h_0$, then $\tilde{E}$ may or not be  irreducible as a complex representation of $\widetilde{\h_0}$.  If it is, we denote by $\mu_E$ its highest weight.  If it is not, then $E$ admits a complex structure as a module over $\h_0$. We then set $\mu_E$ to be the highest weight of $E$ as a complex representation of $\widetilde{\h_0}$.
Recall that the latter case never occurs when $E$ is taken to be $\g_{-1}$.

We are now in position to present the main result.
\begin{thm}\label{thm.casival}
  Let $\g$ be a finite-dimensional graded  maximal  Lie subalgebra of polynomial vector fields. Assume moreover that $\widetilde{\h_0}$ is simple.  Then the Casimir operator $\casi{\Lt}{\g,\beta}$ is diagonalizable. 

Indeed, for all $k\in\N$, its restriction to the sections of an irreducible submodule $I_{k,p}$
of $\vee^k \g_{-1}$ over $\h_0$ equals
\begin{multline}\label{eq.casival}
\frac{1}{2d}(d\delta -k)(d(\delta -1)-k)\\
\mbox{}%
+\frac{\dim(\h_0)}{2(\mu_{\g_{-1}},\mu_{\g_{-1}}+S)d+\dim(\h_0)}(\mu_{I_{k,p}}, \mu_{I_{k,p}}+S)
\end{multline}
times the identity on $I_{k,p}$.
\end{thm}
\begin{proof}
Let $I_{k,p}$ be an irreducible submodule of $\vee^k \g_{-1}$ over $\h_0$.

Assume first that $\K=\C$. Restricted to $I_{k,p}$, 
\begin{equation}\label{eq.casiscal}
\casi{\rhot{k}}{\h_0,\beta_0}= (\mu_{I_{k,p}}, \mu_{I_{k,p}}+S)
\end{equation}
times the identity (see for instance~\cite[p.~122]{hum}). 
As it was recalled in~Proposition~\ref{prop.dualbases}, $\beta\restr{\h_0}$ is nondegenerate.
Under the assumption that $\h_0$ be simple, there exists $l\in\C$ such that 
\begin{equation}
  \label{eq.casi2}
\beta\restr{\h_0}=l\beta_0.  
\end{equation}
 Therefore,
\[
\casi{\rhot{k}}{\h_0,\beta\restr{\h_0}}=\frac{1}{l}\casi{\rhot{k}}{\h_0,\beta_0}.
\]
In order to compute $l$, we recall that $\g_{-1}$ and $\g_1$ are dual $\h_0$-modules and remark that
 for all $x,y\in\h_0$,
\begin{eqnarray*}
   \beta\restr{\h_0}(x,y) &=& 2\tr(\ad(x)\restr{\g_{-1}}\ad(y)\restr{\g_{-1}})+\tr(\ad(x)\restr{\h_0}\ad(y)\restr{\h_0}) \\
   &=& 2 \beta_{\rhot{1}}(x,y)+\beta_0(x,y),
 \end{eqnarray*}
where $\beta_{\rhot{1}}$ is the bilinear form associated to the representation of $\h_0$ on $\g_{-1}$. 

The latter formula also writes 
\[
\beta_{\rhot{1}}=\frac{l-1}{2}\beta_0
\]
and it implies
\[
\casi{\rhot{1}}{\h_0,\beta_0}=\frac{l-1}{2} \casi{\rhot{1}}{\h_0,\beta_{\rhot{1}}}.
\]
Since $\rhot{1}$ is faithful,
\begin{equation}\label{eq.casi3}
\casi{\rhot{1}}{\h_0,\beta_{\rhot{1}}}=\frac{\dim(\h_0)}{d}
\end{equation}
times the identity on $\g_{-1}$.  But, in the same spirit as before,
\[
\casi{\rhot{1}}{\h_0,\beta_0}=(\mu_{\g_{-1}},\mu_{\g_{-1}}+S).
\]
Hence the result over the field of complex numbers.

Let us now handle the case $\K=\R$.  We first adapt Formula~(\ref{eq.casiscal}).        
 The Killing form $\beta_0^\C$ of $\h_0^\C$ restricts to $\beta_0$ on $\h_0\cong\h_0\otimes 1$.  We then look at the Casimir operator associated to the complexification of $\rhot{k}$, $\h_0^\C$ and $\beta_0^\C$.  %This operator is the complex linear extension of $\casi{\rho^k}{\h_0,\beta_0}$.

Two cases may occur.  On the one hand, if $I_{k,p}^\C$ is simple as a complex representation of $\h_0^\C$ --- the latter Casimir operator equals a multiple of the identity, which is real since it is given as above by a scalar product.  It can be computed with the use of the bases of $\h_0$ described in Proposition~\ref{prop.dualbases} and it is thus known on $I_{k,p}$ as well.  On the other hand, if $I_{k,p}^\C$ is reducible, $I_{k,p}$ admits a complex structure.  It becomes a simple complex representation of $\h_0$ and therefore of $\h_0^\C$.

Furthermore, Formulas~(\ref{eq.casi2}) --- with $l\in\R$, this time --- and~(\ref{eq.casi3}) still hold.  Indeed, neither $\h_0$ nor $\rhot{1}$ admit a complex structure as modules over $\h_0$.

Hence the proof.
\end{proof}
The eigenvalue formula~(\ref{eq.casival}) is easily shown to coincide, when $\g$ is taken to be $\alg{sl}{n+1,\R}$, with the formula given in~\cite[Prop. 2]{lecras}.
\section{Building equivariant quantizations}\label{sec.buildeqq}
Throughout this Section, we assume that $\widetilde{\h_0}$ is simple, in order to apply Theorem~\ref{thm.casival}.

For the sake of simplicity, we shall drop the subscripts when writing Casimir operators.  From now on, they will indeed be computed with respect to $\g$ and $\beta$.  We shall furthermore denote by $\casinil$ the difference $\casi{\Lop}{}-\casi{\Lt}{}$.

For all $k\in\N$, let us also decompose 
\begin{equation}\label{eq.decomp}
\vee^k \g_{-1} = \bigoplus_{(p)} I_{k,p}
\end{equation}
into irreducible submodules over $\H_0$, set $E_{k,p}=\Gamma(I_{k,p})$ and denote by $\alpha_{k,p}$ the eigenvalue of $\casi{\Lt}{}$ on $E_{k,p}$.
% We also group the homogeneous eigenvectors by letting
% \[
% \left\{
% \begin{array}{l}
% \widehat{I_{k,p}}\\[1ex]
% \widehat{E_{k,p}}
% \end{array}
% \right.
% =\bigoplus_{q:\alpha_{k,q}=\alpha_{k,p}} 
% \left\{
% \begin{array}{l}
% I_{k,q}\\[1ex]
% E_{k,q}
% \end{array}
% \right.
% \]
% for all $p$ such that $\alpha_{k,p}$ is an eigenvalue of $\casi{\Lt}{}$.

%%
\subsection{The tree-like subspace associated to $\gamma$}
For all $X\in\g_1$, the difference $\gamma(X)=\Lop_X-\Lt_X$ is a differential operator of order~$0$ on $\S$ with constant coefficients.  Identifying tensors in $\vee \g_{-1}$ with tensor fields with constant coefficients, we can thus consider that $\gamma(X)$ is defined on tensors as well.

\begin{lem}\label{lem.subm}
  Let $k\in \N$ and $F$ be a submodule of $\vee^k \g_{-1}$ over $\h_0$.
  Then $\gamma(\g_1)(F)$ is a submodule of $\vee^{k-1} \g_{-1}$ over $\h_0$.
\end{lem}
\begin{proof}
  Since $\gamma$ is a Chevalley-Eilenberg cocycle,  we get the following relation
\[
\Lt_Y\circ\gamma(X)+\gamma(Y)\circ\Lt_X-\Lt_X\circ\gamma(Y)-\gamma(X)\circ\Lt_Y=\gamma([Y,X])
\]
for all vector fields $X$ and $Y$.  If $Y\in\g_0$ and $X\in\g_1$, we deduce
\[
\Lt_Y\circ\gamma(X)P=\gamma(X)\circ\Lt_Y P+\gamma([Y,X])P \in \gamma(\g_1) F
\]
for all $P\in F$.  Hence the proof.
\end{proof}
%%%
%%%
We define the \emph{tree-like subspace associated to $\gamma$}, starting at an irreducible submodule $I_{k,p}$~: 
\[
\tree(I_{k,p})=\bigoplus_{l\in\N}\tree^l(I_{k,p}),
\]
 where $\tree^0(I_{k,p})=I_{k,p}$ and $\tree^{l+1}(I_{k,p})=\gamma(\g_1)(\tree^l(I_{k,p}))$ for all $l\in\N$.  
In view of Lemma~\ref{lem.subm}, these spaces are $\h_0$-submodules. Using Theorem~\ref{thm.casival}, one can compute the spectrum of the restriction of $\casi{\Lt}{}$ to these.
The spaces $\tree^l(E_{k,p})$ are defined in the same way.

Recall that the module structure defined by $\Lop$ is related to two parameters $\lambda$ and $\mu$ and that their difference $\delta$ is called \emph{shift}.
As one expects it to be, the possibility of building equivariant quantization depends on the values of $\lambda$ and $\mu$.  We shall say that an ordered pair of parameters $(\lambda,\mu)$ is \emph{critical} if there exists $E_{k,p}$ such that the corresponding eigenvalue $\alpha_{k,p}$ belongs to the spectrum of the restriction of $\casi{\Lt}{}$ to $\oplus_{l> 0} \tree^l(E_{k,p})$.  In the same fashion, we shall say that a shift value $\delta$ is critical if there exists a value of $\lambda$ such that $(\lambda, \mu)$ is critical in the previous sense.

The following straightforward lemma shows the link between the existence of a $\g$-equivariant quantization and the last definition.
\begin{lem}\label{lem.nc}
  Let $I_{k,p}$ be an irreducible submodule of $\vee^k \g_{-1}$ and $\g_{-1}^*\otimes I_{k,p}$ denote the subspace of $E_{k,p}$ made up of sections with linear coefficients.  Then
\[
\casinil(\g_{-1}^*\otimes I_{k,p})=\gamma(\g_1)(I_{k,p}).
\]
\end{lem}
\begin{proof}
In the basis $(e_i)$ of $\g_{-1}$ chosen in Proposition~\ref{prop.dualbases}, $\Lt_{e_i}$ takes the local form $\del{i}$.
  It then follows that
\[
\casinil(\sum_{l,m} a^l_m x^m u_l)=2\sum_{l,m} a^l_m\gamma(\epsilon^m)u_l,
\]
for all $a_m^l\in\K$ and $u_l\in I_{k,p}$, $(1\leq l,m\leq d)$.
\end{proof}
In the same fashion, the following lemma is immediate.
\begin{lem}
  For all $u\in E_{k,p}$, $\casinil(u)\in\gamma(\g_1)(E_{k,p})$.
\end{lem}
\begin{thm}\label{thm.q}
  If $(\lambda,\mu)$ is not critical, then there exists a $\g$-equivariant quantization.
\end{thm}
\begin{proof}
 The proof machinery goes as in~\cite{dlo}.  We give it for the sake of completeness.

Let $P\in E_{k,p}$.  We first prove that there exists a unique $\hat{P}\in\tree(E_{k,p})$ such that $P$ is the principal symbol of $\hat{P}$ and that $\hat P$ is an eigenvector of $\casi{\Lop}{}$ associated to the eigenvalue $\alpha_{k,p}$.
For all $R\in\S$, write $R_l$ the projection of $R$ onto $\S^l$.
Now, for all $l<k$, $\tree^{k-l}(E_{k,p})$ splits into a sum of eigenspaces $W_{l,q}$ associated to the eigenvalues $\beta_{l,q}$ of $\casi{\Lt}{}$. Denote by $R_{l,q}$ the projection of $R$ onto $W_{l,q}$.

With these notations, the equation $\casi{\Lop}{}\hat{P}=\alpha_{k,p} \hat{P}$ can be rewritten :
\[
\left\{
  \begin{array}{l}
    \casi{\Lt}{} P=\alpha_{k,p} P\\
    \sum_q \left( (\casinil(\hat{P}_{l+1}))_{l,q} + \casi{\Lt}{} \hat{P}_{l,q}\right) = \alpha_{k,p} \hat{P}_{l,q},\quad \forall l,q,
  \end{array}
\right.
\]
which reduces to
\begin{equation}\label{eq.quantsyst}
\left\{
  \begin{array}{l}
    \casi{\Lt}{} P=\alpha_{k,p} P\\
    (\alpha_{k,p}-\beta_{l,q}) \hat{P}_{l,q}=(\casinil(\hat{P}_{l+1}))_{l,q}, \quad \forall l,q,
  \end{array}
\right.
\end{equation}
where the last equation must be satisfied for all $l<k$ and $q$.
The existence and the properties of the correspondence $P\mapsto\hat P$ follow from the observation that the latter system is triangular and admits a unique solution.  Indeed, the right hand side of the equations involving $\casinil$ always belongs to $\tree(E_{k,p})$.

Now, let $\mathcal{Q}$ denote the linear extension of this correspondence.  It remains to prove that it is equivariant with respect to $\g$.
 It suffices to check that 
\[
\Lop_X\circ \mathcal{Q}(P)=\mathcal{Q}\circ\Lt_X(P),
\]
for all $X\in\g$, all $k\in\N$ and  all eigenvectors $P\in\S^k$ of $\casi{\Lt}{}$ associated to any eigenvalue $\alpha_{k,p}$.  But both sides of this condition are eigenvectors of $\casi{\Lop}{}$ associated to the same eigenvalue $\alpha_{k,p}$.  Moreover, they have the same principal symbol: $\Lt_X(P)$.  Since $E_{k,p}$ and $\tree(E_{k,p})$ are respectively closed under $\Lt_X$ and $\Lop_X$, both sides belong to the latter tree.  The first part of the proof ensures that they coincide.
\end{proof}

\section{Examples}\label{sec.examples}
We shall now apply the method described in the previous section to two particular algebras.  The treatment will be done in a concurrent way. 
Throughout this Section, $\g$ denotes one of the algebras $\O$ and $\Sp$ defined below.

 We shall significantly refine Theorem~\ref{thm.q} by proving that $0$ is never a critical shift value and obtaining the uniqueness of the quantization.

\subsection{Orthogonal and Symplectic algebras}
%The algebras we shall consider are tangent to the action of the linear symplectic (resp. orthogonal) group on the Lagrangian (resp. orthogonal) Grassmannian in a vector space endowed with a symplectic form (resp. scalar product).
From now on, we assume that $n$ is an integer greater than $1$.
It is well known that the Lie subalgebras $\alg{so}{n,n,\K}$ and $\alg{sp}{2n,\K}$ of the general linear algebra $\alg{gl}{2n,\K}$ can be realized as $3$-graded algebras.  These are described for instance in~\cite[pp. 893-894]{koba1}.

For the constructions below to be self-contained, we only need to recall that $\alg{so}{n,n,\K}$ is written 
\[
\O=\O_{-1}\oplus\O_0\oplus\O_1,
\]
where $\O_{-1}=\wedge^2\K^n$, $\O_1=\wedge^2\K^{n*}$ and $\O_0=\alg{gl}{n,\K}$.  
For all $A\in\O_0$ and $h\in\O_{-1}\oplus\O_{1}$, 
\[
[A,h]=\rho(A)h,
\]
where $\rho$ is the natural representation of $\O_0$ on $\O_{-1}\oplus\O_1$.
The Euler element is the identity transformation of $\O_{-1}$. 
We shall refer to $\O$ as the \emph{orthogonal} algebra.

Similarly, $\alg{sp}{2n,\K}$ is written 
\[
\Sp=\Sp_{-1}\oplus\Sp_0\oplus\Sp_1,
\]
where $\Sp_{-1}=\vee^2\K^n$, $\Sp_1=\vee^2\K^{n*}$ and $\Sp_0=\alg{gl}{n,\K}$. 
The same statements about the bracket and Euler element hold.  We shall refer to this algebra as the \emph{symplectic} algebra.

\subsection{Casimir operator eigenvalues}
In the examples under consideration, the subalgebra $\h_0$ is isomorphic to $\alg{sl}{n,\K}$ and $\widetilde{\h_0}=\alg{sl}{n,\C}$ is obviously simple. The data introduced to state Theorem~\ref{thm.casival} are classical.
Let us denote by  $\alg{d}{n,\K}$ the matrix subalgebra of diagonal matrices of $\alg{gl}{n,\K}$ and $D_j\in\alg{d}{n,\K}, (j=1,\ldots,n-1)$, as the diagonal matrix 
\[
\mathrm{diag}(0,\ldots,0,\overbrace{1}^{j},0,\ldots,0,-1).
\]
These diagonal matrices generate the Cartan subalgebra $\alg{sl}{n,\K}\cap \alg{d}{n,\K}$ of $\alg{sl}{n,\K}$. In its dual space, we define $\delta_j$ by $\delta_j(D_i)=\delta_{ij}$ for all $i,j\in\{1,\ldots,n-1\}$.  As it is common, we set $\delta_n=-\sum_{i=1}^{n-1}\delta_i$ as well.  Then, a simple root system of $\alg{sl}{n,\K}$ is given by the $\delta_i-\delta_{i+1}$, $(i=1,\ldots,n-1)$.
The sum $S$ of the positive roots with respect to this system equals $2\sum_i (n-i)\delta_i$.
The Killing form $\beta_0$ of $\alg{sl}{n,\K}$ is given by $\beta_0(A,B)=2n\,\tr(AB)$ for all $A,B\in\alg{sl}{n,\K}$.  The induced scalar product satisfies
\begin{equation}\label{eq.scalprod}
(\delta_i,\delta_j)=\frac{1}{2n^2}(n\delta_{ij}-1)\text{ and }(\delta_i,S)=\frac{n-2k+1}{2n},
\end{equation}
for all $i=1,\ldots,n$.

\begin{figure}[htbp]
  \begin{center}
    \scalebox{.5}{ \includegraphics{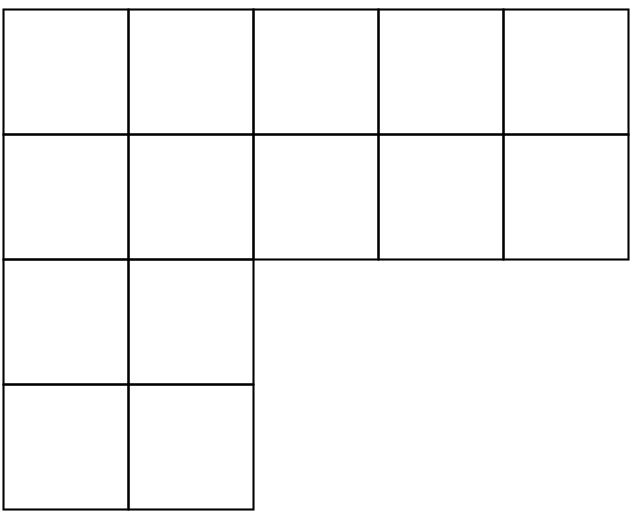}}
    \quad
    \scalebox{.5}{\includegraphics{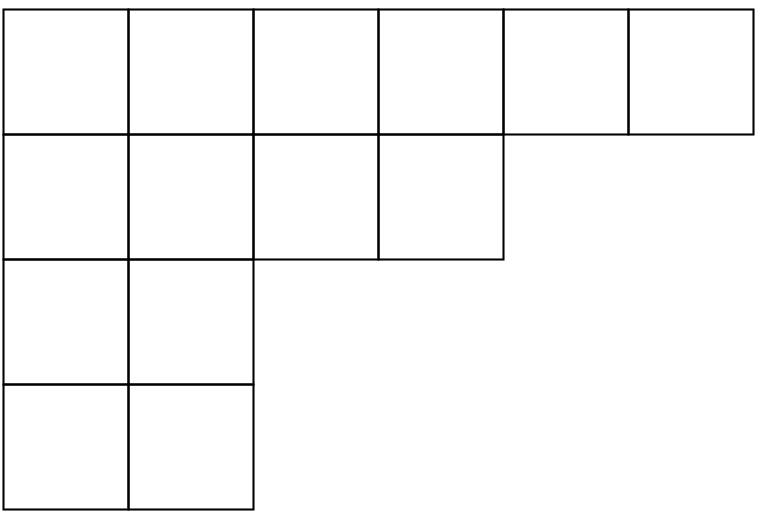}}
    
  \end{center}
  \caption{Ferrers diagrams of irreducible submodules of $\vee^k\O_{-1}$ and $\vee^k\Sp_{-1}$ $(n\geq4)$.}
\label{fig.Ferrers}
\end{figure}

Now, let $\K=\C$.  The decomposition of $\vee^k\g_{-1}$ into irreducible submodules over $\g_0$ is given in~\cite{gw}.  Those submodules may be generated by the action of real matrices on some (real) highest weight vector.  Therefore, those results can be used when $\K=\R$ as well.
As it is well-known, irreducible submodules can be conveniently indexed by Ferrers diagrams, which in turn can be denoted by elements of $\N^n$.
We shall respectively denote by  $(5,5,2,2)$ and $(6,4,2,2)$  the Ferrers diagram given in Figure~\ref{fig.Ferrers}. They admit $5(\delta_1+\delta_2)+2(\delta_3+\delta_4)$ (resp. $6\delta_1+4\delta_2+2(\delta_3+\delta_4)$) as highest weight.  

Theorem~5.2.11  in~\cite{gw} states that 
$\vee^k\O_{-1}$ splits as a sum of one copy of each  irreducible submodule of weight  $\mu=\sum_{i=1}^nk_i\delta_i$ such that 
\begin{enumerate}
 \item $k_1\geq\cdots\geq k_n\geq 0$, ($\mu$ is dominant nonnegative),
 \item  $\sum_i k_i=2k$,
 \item $k_{2i-1}=k_{2i},\quad\forall i\in\{1,\ldots,\lfloor n/2\rfloor\}$,
 \item $k_{n}=0$ if $n$ is odd.
\end{enumerate}
Theorem~5.2.9 ibidem states that 
$\vee^k\Sp_{-1}$ splits as a sum of one copy of each irreducible submodule of weight 
 $\mu=\sum_{i=1}^n k_i\delta_i$ such that 
\begin{enumerate}
 \item $k_1\geq\cdots\geq k_n\geq 0$, ($\mu$ is dominant nonnegative),
 \item  $\sum_i k_i=2k$,
 \item $k_i\in 2\N,\quad\forall i\in\{1,\ldots,n-1\}$.
\end{enumerate}

Let us compute explicitly the value of Formula~(\ref{eq.casival}).
For all submodules $R$ of $\O$ or $\Sp$ with highest weight $\mu_R$ described by a Ferrers diagram~$(k_1,\ldots,k_n)$, Formula~(\ref{eq.scalprod}) shows that
\begin{equation}\label{eq.scalprod2}
(\mu_R,\mu_R+S)=\frac{1}{2n^2}%\left(
\sum_{i,j=1}^n ( k_i k_j (n\delta_{ij}-1)+2 k_i(n-j)(n\delta_{ij}-1)).
%\right).
\end{equation}

In the \emph{orthogonal} case, $d=\frac{n(n-1)}{2}$ and the highest weight of $\g_{-1}=\O_{-1}$ is $\delta_1+\delta_2$.  Let $R$ now denote an irreducible submodule of $\vee^k\O_{-1}$ associated to a Ferrers diagram $\vec{k}=(k_1,\ldots,k_n)$.  A direct computation using~(\ref{eq.casival}) and~(\ref{eq.scalprod2}) shows that the eigenvalue of $\casi{\Lt}{}$ associated to $R\otimes\Delta^\delta\O_{-1}$ equals
\begin{equation}\label{eq.orthoval}
\alpha_{o}(\vec k)=\frac{n(n-1)}{4}\delta^2-(k+\frac{n(n-1)}{4})\delta+\frac{n}{n-1}k+\frac{\sum_{i=1}^n k_i(k_i-2i)}{4(n-1)}.
\end{equation}

In the \emph{symplectic} case, $d=\frac{n(n+1)}{2}$ and the highest weight of $\g_{-1}=\Sp_{-1}$ is $2\delta_1$.  Let $R$ now denote an irreducible submodule of $\vee^k\Sp_{-1}$ associated to a Ferrers diagram $(k_1,\ldots,k_n)$.   Then the eigenvalue of $\casi{\Lt}{}$ associated to $R\otimes\Delta^\delta\O_{-1}$ equals
\begin{equation}\label{eq.symplecval}
\alpha_{s}(\vec k)=\frac{n(n+1)}{4}\delta^2-(k+\frac{n(n+1)}{4})\delta+ k +\frac{\sum_{i=1}^n k_i(k_i-2i)}{4(n+1)}.
\end{equation}

\subsection{Another tree}
In both symplectic and orthogonal cases, it is easy to check that the difference of two eigenvalues corresponding to different degrees $k$ cannot be identically zero.  Indeed, such a difference is a linear expression in $\delta$ with rational coefficients.  Thus there exist infinitely many values of the shift for which a quantization exists.

We shall now develop two important refinements. First, we shall exclude $0$ from the critical shift values in both the symplectic and orthogonal cases.  Then, given any non critical value of the shift, we shall prove that the only equivariant quantization is the one we have built.

In order to prove that $0$ is not critical, it is unfortunately not sufficient to check all the eigenvalues by a straight inspection.  For instance,  the eigenvalues associated to diagrams $(6,2,2,2)$ and $(6,4)$ are equal when $n=6$ (resp. $n=5$) in the orthogonal (resp. symplectic) case.
But it is clear from Equation~(\ref{eq.quantsyst}) that only some of those equalities can actually prevent the quantization from existing.
% Let us give a simple remark.
% Recall that $\h_0$ is made up of linear vector fields over $\g_{-1}$ while $\H_0$ is the corresponding algebra of linear transformations of that space.  
% Let $\rho$ denote the natural representation of $\H_0$ on $\g_{-1}$ and $\rho^*$ its natural extension to $\g_{-1}^*$.
% Let also $W_0$ be the canonical fiber of a trivial vector bundle $W$ over $\g_{-1}$ and let it be endowed with a representation $\rho_{W_0}$ of $\H_0$.  Then the representation defined by the Lie derivative with respect to (linear) vector fields in $\h_0$ on $\g_{-1}^*\otimes W_0$ --- viewed as the space of sections with linear coefficients of $W$ --- is trivially isomorphic to the representation $\rho^*\otimes\rho_{W_0}$ of $\H_0$ on $\g_{-1}^*\otimes W_0$. 

Let $I\subset\vee^k \g_{-1}$ be an irreducible submodule over $\H_0$.
We define a bigger tree than $\tree(I)$ as follows.
Let $\treebig^1(I)$ be the sum of all irreducible submodules $J_p$ in $\vee^{k-1} \g_{-1}$ that are isomorphic to an irreducible submodule of $\g_{-1}^*\otimes I$. Define $\treebig^2(I)=\bigoplus_{(p)} \treebig^1(J_p)$ and continue recursively.  We write
\[
\treebig(I)=I\oplus\bigoplus_{k\geq 1}\treebig^k(I).
\]
Consider now the natural representation of $\H_0$ on $\g_{-1}^*\otimes I$.  It is isomorphic to the representation defined by
the Lie derivative in the direction of (linear) vector fields of $\h_0$  on the space of sections valued in $I$ with linear coefficients (cf. Subsection~\ref{subsec.diag}).
Lemma~\ref{lem.nc} and the invariance of $\casinil$ under $\h_0$ then allow to conclude that for all $\lambda$, $\tree(I)$ is indeed a subset of $\treebig(I)$. 

It is customary to order  the Ferrers diagrams as follows:
\[
\vec k\leq \vec l \equiv (k_i\leq l_i, \forall i\leq n)
\]
and of course
\[
\vec k < \vec l\equiv (\vec k\leq \vec l \text{ and } \vec k\neq \vec l).
\]

Then we can describe $\treebig(I)$ in the examples under consideration.
\begin{lem}\label{lem.lr}
  Let $K\subset\vee^k\g_{-1}$ be an irreducible submodule over $\H_0$ whose type is given by the Ferrers diagram $\vec k$.  Then an irreducible submodule $L\subset\vee^l\g_{-1}$ with type $\vec l$, $(l<k)$ is in $\treebig(K)$ only if $\vec l< \vec k$.
\end{lem}
\begin{proof}
  It suffices to determine the diagrams occurring in the decomposition
  of $\g_{-1}^*\otimes K$ into irreducible components using
  Littlewood-Richardson rule (see for instance~\cite[pp. 455-456]{lr}). 

Let us explicit the proof in the orthogonal case, for which $\g_{-1}^*$ is represented by a column of height $n-2$ and width $1$.  The irreducible components of $K\otimes \g_{-1}^*$ are then associated to diagrams made up by adding one box to $n-2$ rows of the diagram associated to $K$.

Then, one needs to know which of these new diagrams represent irreducible components isomorphic to one occurring in the decomposition of $\vee^{k-1}\g_{-1}$.  But the later admit diagrams with $2(k-1)$ boxes while the former have $2k+n-2$.  In order to describe isomorphic $\alg{sl}{n,\K}$ submodules, they should differ by one column of height $n$ and width $1$ on the left. The diagram with $2k-n+2$ boxes may thus only be isomorphic to  a diagram smaller than the original diagram of $K$.  The conclusion follows by induction.
\end{proof}

\begin{thm}\label{thm.notzero}
All critical shift values belong to the set
\[
\{\frac{n}{n-1}+\frac{\sum_{i=1}^n (k_i-l_i)(k_i+l_i-2i)}{4(n-1)(k-l)}:\vec k>\vec l\}
\]
in the orthogonal case and
\[
\{ 1+\frac{\sum_{i=1}^n (k_i-l_i)(k_i+l_i-2i)}{4(n+1)(k-l)}:\vec k>\vec l\}
\]
in the symplectic case, where $\vec k$ and $\vec l$ describe admissible Ferrers diagrams.
They are greater than $0$.  In particular, there exists a $\g$-equivariant quantization into operators that preserve the weight of their arguments.
\end{thm}
\begin{proof}
    Assume that two eigenvalues associated to $K$ and $L$ taken as above are equal.  Then the shift value is fixed to a positive value.
  
For instance, in the orthogonal case, the announced value is 
 not less than
\[
\frac{n}{n-1}+\frac{\sum_i (k_i^2-l_i^2)}{4(n-1)(k-l)}-2n \frac{\sum_i (k_i-l_i)}{4(n-1)(k-l)},
\]
 which is greater than $0$.  Indeed, the last term sums up to the first and there exists an index $i$ such that $k_i>l_i$.

The proof goes the same way in the symplectic case.  Hence the conclusion.
\end{proof}

Let us now turn to the uniqueness problem. 
\begin{lem}\label{lem.unique}
  Assume that $\delta$ does not belong to a set presented in Theorem~\ref{thm.notzero}. Let $k,l\in\N$ such that $l<k$.  Then there exists no $\g$-equivariant map from  $(\S^k,\Lt)$ to $(\S^l,\Lt)$.
\end{lem}
\begin{proof}
  Assume that $T$ is such a map.  Roughly speaking, we shall obtain a contradiction about the value of  $\delta$ because $T$ must preserve the tree-like subspaces $\treebig$.

As it was proved in~\cite[Lemma~7.1]{lecsln}, the equivariance of $T$ with respect to $\g_{-1}$ (i.e. every constant vector field) and $\euler$ implies that it is a differential operator with constant coefficients.  We can thus write
\[
T=\sum_{r=0}^R T_r,
\]
with $T_r$ an homogeneous differential operator of order $r$.

In view of~(\ref{eq.lieeuler}), $[\Lt_\euler, T]=0$ leads furthermore to
\[
\sum_{r=0}^R (k-l-r) T_r=0
\]
and, therefore, $T=T_{k-l}$.

Let now, as in~(\ref{eq.decomp}), $I_{k,p}$ be an irreducible submodule of $\vee^k\g_{-1}$ over $\H_0$, described by $\vec{\imath}$.  We know that $T$ is entirely defined by its values on the sections in $\Gamma(I_{k,p}\otimes \Delta^\delta \g_{-1})$ with polynomial coefficients of degree $k-l$. We recall that the Lie derivative in the direction of vector fields in $\h_0$ on those has no effect on the ``density part'' and corresponds to the natural representation of $\H_0$ on 
\begin{equation}\label{eq.polysec}
\vee^{k-l} \g_{-1}^*\otimes I_{k,p}.
\end{equation}
The image of such sections through the application of $T$ is made of sections with constant coefficients.  This image corresponds to a submodule $F$ of $\vee^l(\g_{-1})$ over $\H_0$.
The irreducible components of $F$ necessarily appear in the decomposition of~(\ref{eq.polysec}) and thus in that of $\otimes^{k-l}\g_{-1}^*\otimes I_{k,p}$.  

Our last argument goes as in the proof of Lemma~\ref{lem.lr}.  Let $\vec f$ describe a submodule of $F$ isomorphic to a submodule $L$ of $\otimes^{k-l}\g_{-1}^*\otimes I_{k,p}$. Let $\vec l$ be the diagram describing $L$. On the one hand, $\vec f$ is obtained by removing $2(k-l)$ (resp. $(k-l)$) columns on the left of $\vec l$ in the symplectic (resp. orthogonal) case.  On the other hand, in application of the Littlewood-Richardson rule, $\vec l$ is obtained by adding $2(k-l)(n-1)$ (resp. $(k-l)(n-2)$) boxes  to $\vec{\imath}$, with no more than $2(k-l)$ (resp. $(k-l)$) boxes in a single row.  

Therefore, $\vec f < \vec{\imath}$. But the invariance of $T$ ensures that the values of $\casi{\Lt}{}$ on $F$ and $I_{k,p}$ coincide, which contradicts the hypothesis on $\delta$.
\end{proof}
\begin{cor}
  If the shift is chosen as in Theorem~\ref{lem.unique} then the $\g$-equivariant quantization is unique.
\end{cor}
\section*{Acknowledgments}
We are very grateful to  M.~De~Wilde, C.~Duval, P.~Lecomte  and V.~Ovsienko for their 
interest and suggestions.
The first author should like to thank W. Bertram for interesting discussions and the  Belgian 
National Fund for Scientific Research (FNRS) for his Research 
Fellowship.

\bibliographystyle{plain}
\bibliography{QuantumCasimirs}
%\nocite{*}

\noindent Institute of Mathematics, B37
\\University of Li\`ege,
\\B-4000 Sart Tilman
\\Belgium
\\
\\Email addressses:
\\f.boniver@ulg.ac.be
\\p.mathonet@ulg.ac.be
\end{document}